\documentclass{amsart}
\usepackage{amssymb,amsmath,latexsym}
\begin{document}

\title{\bf A remark on the Chebotarev theorem about roots of unity.}
\author{F. Pakovich}
\thanks{Research supported by the ISF, Grant No. 979/05}

\date{}
\address{Department of Mathematics, Ben Gurion University,
P.O.B. 653, Beer-Sheva 84105, Israel}

\email{pakovich@math.bgu.ac.il}

\subjclass{11T22, 11C08}
\keywords{Chebotarev theorem, roots of unity, polynomials}

\begin{abstract}
Let $\Omega$ be a matrix 
with entries 
$a_{i,j}=\omega^{ij},$ $1\leq i,j \leq n,$ where $\omega=e^{2\pi \sqrt{-1}/n},$ $n\in \mathbb N.$ 
The Chebotarev theorem states that if $n$ 
is a prime then any minor of $\Omega$ is non-zero. 
In this note we provide an analogue of this statement for composite $n.$

\end{abstract}

\maketitle

\def\d{{\rm d}}
\def\D{{\rm D}}
\def\I{{\rm I}}

\def\C{{\mathbb C}}
\def\N{{\mathbb N}}
\def\P{{\mathbb P}}   
\def\Z{{\mathbb Z}}
\def\d{{\rm d\,}}
\def\deg{{\rm deg\,}}
\def\Det{{\rm Det}}\def\dim{{\rm dim\,}}
\def\Ker{{\rm Ker\,}}
\def\Gal{{\rm Gal\,}}
\def\St{{\rm St\,}}
\def\Sym{{\rm Sym\,}}
\def\Mon{{\rm Mon\,}}

Let $\Omega$ be a matrix 
with entries 
$a_{i,j}=\omega^{ij},$ $1\leq i,j \leq n,$ where $\omega=e^{2\pi \sqrt{-1}/n},$ $n\in \mathbb N.$ 
The Chebotarev theorem states that if $n$ 
is a prime then any minor of $\Omega$ is non-zero. Chebotarev's proof 
of this theorem and the references to other proofs can be found in \cite{l}. Yet another 
proofs can be found in recent papers \cite{g}, \cite{t}.

For a complex polynomial $P(z)$ denote by $w(P)$ the number of non-zero coefficients of $P(z).$
It is easy to see that the Chebotarev theorem is equivalent to the following statement: if a non-zero polynomial $P(z),$ $\deg P(z) \leq n-1,$  
has $k$ different roots which are $n$-roots of unity 
then $w(P)> k$ whenever $n$ is a prime. 

A natural question is: how small can be $w(P)$ if $n$ is a composite number ?
The example of the polynomial $$D_{n,r,l}(z)=z^l(1+z^r+z^{2r}+ ... + z^{(\frac{n}{r}-1)r}),$$ where $r\vert n$, $0\leq l \leq r-1,$ shows that $w(P)$ could be as small as $n/(n-k).$  
In this note we show that actually it is the ``worst'' possible case.

\vskip 0.2cm
\noindent{\bf Theorem.} {\it Let $n$ be a composite number and $P(z)$ be a non-zero complex polynomial, $\deg P(z) \leq n-1.$ 
Suppose that $P(z)$ has
exactly   
$k$ different roots which are $n$-roots of unity. Then the inequality 
$$ w(P)\geq  \frac{n}{n-k}\eqno(*)$$ holds. Furthermore, the equality attains if and only 
if $P(z)$ up to a multiplication by a complex number coincides with
$D_{n,r,l}(\omega^j z)$ 
%$$P(z)=\alpha z^l(1+(\nu^j)z^r+ (\nu^j)^{2}z^{2r}+ ... + (\nu^j)^{(\frac{n}{r}-1)}z^{(\frac{n}{r}-1)r}),$$ 
for some $j,$ $0\leq j \leq n-1,$ and $r,l$ as above.  
}
\vskip 0.2cm
\noindent{\it Proof.} Let $P(z)=p_{0}+p_{1}z+...+p_{n-1}z^{n-1}$ and let 
$$C=\begin{pmatrix}p_0 & p_1 & ... & p_{n-1}\\
p_{n-1} & p_0 & ...& p_{n-2}\\
... & ... & ... & ...\\
p_1& p_2 & ... & p_0\end{pmatrix}$$ be the circulant matrix generated by the coefficients 
of $P(z).$ We will denote the row vectors of $C$ by $\vec{t_j},$ $0\leq j \leq n-1$.
Set $r={\rm rk}\, C.$
The key observation is that the number $k$ is equal to the number $n-r$. To establish it notice that eigenvectors of $C$ are
$\vec{f_i}=((\omega^i)^0,(\omega^i)^1,...,(\omega^i)^{(n-1)}),$ $0\leq i \leq n-1,$ and the corresponding eigenvalues are $P(\omega^{i}),$ $0\leq i \leq n-1.$ Furthermore, the vectors $\vec{f_i},$ $0\leq i \leq n-1,$ form a basis of $\mathbb C^n$. The matrix $C$ is diagonal with respect to this basis and therefore 
$k=n-r$.

It follows that in order to prove inequality (*) it is enough to establish
the inequality $$w(P)\hskip 0.05cm r\geq n. \eqno(**)$$ This inequality essentially is a particular case 
of Theorem C in \cite{g} and can be established as follows (\cite{g}).
Let $V$ be a vector space generated by the vectors $\vec{t}_j,$ $0\leq j \leq n-1$, and
$R$ be a subset of 
$\{\vec{t}_0, \vec{t}_1, ... , \vec{t}_{n-1}\}$
consisting of $r$ vectors which generate $V.$ Clearly, 
for any $i,$ $1\leq i \leq n,$ there exists a vector $v\in V$ for which
its $i$-th coordinate is distinct from zero. Since
each vector from $R$ has exactly $w(P)$
non zero coordinates it follows that (**) holds.

For a vector $\vec{v}\in \mathbb C^n$ denote by ${\rm supp}\{\vec{v}\}$ the set consisting of numbers $i,$ 
$1\leq i \leq n,$ for which $i$-th coordinate of $\vec{v}$ is non-zero.
Observe now that the equality in (**) attains only if
for any two vectors $v_1,v_2\in R$ we have ${\rm supp}\{\vec{v_1}\}\,\cap\,{\rm supp}\{\vec{v_2}\}=\emptyset$. This implies easily
that ${\rm supp}\{\vec{t_0}\}$ consists of numbers 
congruent by modulo $r$ to a number $l,$ $0\leq l \leq r-1,$ 
and therefore $P(z)=z^lQ(z^r)$ for some polynomial $Q(z)=q_{0}+q_{1}z+...+q_{n/r-1}z^{n/r-1}$
and number $l,$ $0\leq l \leq r-1.$

Furthermore, since the vectors $t_0,t_r,t_{2r}, ... , t_{n/r-1}$ have equal supports
the equality in (**) implies that any two of them are proportional. Therefore, the rank of the circulant matrix $W$ generated by the coefficients of $Q(z)$ equals 1. This implies that
the vector $\vec{q}=\{q_{0},q_{1},...,q_{n/r-1}\}$ is orthogonal to $n/r-1$ vectors from the   
collection $\vec{g_j}=((\nu^j)^0,(\nu^j)^1,...,(\nu^j)^{(n/r)-1}),$ $0\leq j \leq n/r-1,$ 
where $\nu=\omega^r.$ Since $\vec{g_j},$ $0\leq j \leq n/r-1,$ are linearly independent
this implies that there exists $\alpha\in \mathbb C$ such that $\vec{q}=\alpha \vec{g_j}$ for some $0\leq j \leq n/r-1.$

\bibliographystyle{amsplain}

\end{document}